\documentclass[24pt,reqno]{article}
\usepackage
[colorlinks=true, linkcolor=webgreen, filecolor=webbrown,
citecolor=webgreen]{hyperref}

\usepackage{amsmath}
\usepackage{amssymb}
\usepackage{color}
\usepackage{epsf}

\begin{document}
\makeatletter

\begin{center}
\epsfxsize=10in
\end{center}

\def\endofproofmark{$\Box$}

\def\RR{{\mathbb R}}
\def\NN{{\mathbb N}}

\def\endofproofmark{$\Box$}

\begin{center}
\vskip 1cm {\LARGE\bf Some Refinements of Discrete Jensen's}
\vskip 10mm {\LARGE\bf Inequality and Some of Its Applications}
\vskip 1cm \large J. Rooin
\\
\vskip .5cm
Department of Mathematics\\
Institute for Advanced
Studies in Basic Sciences\\
P.O. Box 45195-159\\
Zanjan, Iran\\
\href{mailto:rooin@iasbs.ac.ir}{\tt rooin@iasbs.ac.ir}\\
\end{center}

\date{}

\newtheorem{theo}{Theorem}[section]
\newtheorem{prop}[theo]{Proposition}
\newtheorem{lemma}[theo]{Lemma}
\newtheorem{cor}[theo]{Corollary}
\newtheorem{problem}{Problem}
\def\frameqed{\framebox(5.2,6.2){}}
\def\deshqed{\dashbox{2.71}(3.5,9.00){}}
\def\ruleqed{\rule{5.25\unitlength}{9.75\unitlength}}
\def\myqed{\rule{8.00\unitlength}{12.00\unitlength}}
\def\qed{\hbox{\hskip 6pt\vrule width 7pt height11pt depth1pt\hskip 3pt}
\bigskip}
\newenvironment{proof}{\trivlist\item[\hskip\labelsep{\bf Proof}:]}{\hfill
 $\frameqed$ \endtrivlist}
\newcommand{\COM}[2]{{#1\choose#2}}

\thispagestyle{empty} \null \addtolength{\textheight}{1cm}

\begin{abstract}

In this paper, using  some aspects of convex functions, we refine
discrete Jensen's inequality via weight functions. Then, using
these results, we give some applications in different abstract
spaces and obtain some new interesting inequalities.

\end{abstract}
\section{Introduction}
Jensen's inequality is sometimes called the king of inequalities
[4] because it implies at once the main part of the other
classical inequalities (e.g. those by H\"{o}lder, Minkowski,
Young, and the AGM inequality, etc.). Therefore, it is worth
studying it thoroughly and refine it from different points of
view. There are numerous refinements of Jensen's inequality, see
e.g. [3-5] and the references in them. In this paper, introducing
suitable weight functions, first we give some refinements of
discrete Jensen's inequality, and then using these refinements, we
give several important applications in
various abstract spaces, which extends the results obtained recently [5,6].\\
Throughout this paper, we suppose that $C$ is a convex subset of a
real vector space, $x_1,\cdots,x_n\in C$, and
$\varphi:C\rightarrow{\Bbb R}$ a convex mapping. Also, we suppose
that $\mu=(\mu_1,\cdots,\mu_m)$ and
$\lambda=(\lambda_1,\cdots,\lambda_n)$ are two probability
measures; i.e. $\mu_i,~\lambda_j\geq 0~(1\leq i \leq m,~1\leq j
\leq n)$ with
$$\sum_{i=1}^m \mu_i=1\hspace{1cm}\mbox{and}\hspace{1cm}\sum_{j=1}^n \lambda_j=1.
$$
By a (discrete separately) weight function (with respect to $\mu$
and $\lambda$), we always mean a mapping
$$
\omega:\{(i,j):~1\leq i \leq m,~1\leq j \leq n\}\rightarrow
[0,\infty),
$$
such that
$$
\sum_{i=1}^m\omega(i,j)\mu_i=1~~~~~~(j=1,\cdots,n),
$$
and
$$
\sum_{j=1}^n\omega(i,j)\lambda_j=1~~~~~~(i=1,\cdots,m).
$$
For example, if $u=(u_1,\cdots,u_m)$ and $v=(v_1,\cdots,v_n)$ with
$\|u\|=(\sum_{i=1}^m u_i^2)^{1/2}\leq 1$ and $\|v\|=(\sum_{j=1}^n
v_j^2)^{1/2} \leq 1$ belong to $\mu^\perp$ and $\lambda^\perp$
respectively, then the function  $\omega$ with
$\omega(i,j)=1+u_iv_j~(1\leq i \leq m,~1\leq j \leq n)$ is a
weight function.\\
Also, we say that a quadratic matrix $A=[a_{ij}]_{n\times n}$ with
nonnegative entries is a double stochastic matrix if the sum of
each of its rows and columns is unit; that is
$$
\sum_{i=1}^n a_{ij}=1~~~~~~(j=1,\cdots,n),
$$
and
$$
\sum_{j=1}^na_{ij}=1~~~~~~(i=1,\cdots,n).
$$
If $\omega_1$ and $\omega_2$ are two weight functions, we denote
by $\phi_{\omega_1,\omega_2}$ the real-valued function
\begin{eqnarray}
\phi_{\omega_1,\omega_2}(t)=\sum_{i=1}^m\mu_i\varphi
\left(\sum_{j=1}^n[(1-t)\omega_1(i,j)+t\omega_2(i,j)]\lambda_jx_j
\right)~~~~~~(0\leq t \leq 1),
\end{eqnarray}
and also, if $B=[b_{ij}]_{n\times n}$ and $C=[c_{ij}]_{n\times n}$
are two double stochastic matrices, we put
\begin{eqnarray}
\phi_{B,C}(t)=\frac{1}{n}\sum_{i=1}^n\varphi
\left(\sum_{j=1}^n[(1-t)b_{ij}+tc_{ij}]x_j \right)~~~~~~(0\leq t
\leq 1).
\end{eqnarray}
The aim of this paper is to refine discrete Jensen's inequality
$$
\varphi\left(\sum_{j=1}^n\lambda_jx_j\right)\leq
\sum_{j=1}^n\lambda_j\varphi(x_j)
$$
via weight functions and give some applications in different
abstract spaces which extend the results of [5] and [6].
\section{Refinements} In this section, using the above notations, we refine
discrete Jensen's inequality by
one and two weight functions which extend the results of [5].\\\\
{\bf Lemma 2.1.}~{\it If $\omega$ is a weight function, then
\begin{eqnarray}
\varphi\left(\sum_{j=1}^n\lambda_jx_j\right)
\leq\sum_{i=1}^m\mu_i\varphi\left(\sum_{j=1}^n\omega(i,j)\lambda_jx_j\right)
\leq\sum_{j=1}^n\lambda_j\varphi(x_j).
\end{eqnarray}
In particular, if $A=[a_{ij}]_{n\times n}$ is a double stochastic
matrix, we have
\begin{eqnarray}
\varphi\left(\frac{x_1+\cdots+x_n}{n}\right)
\leq\frac{1}{n}\sum_{i=1}^n\varphi\left(\sum_{j=1}^na_{ij}x_j\right)
\leq\frac{\varphi(x_1)+\cdots+\varphi(x_n)}{n}.
\end{eqnarray}}\\
{\it Proof}.~By the convexity of $\varphi$, we have
$$
\sum_{i=1}^m\mu_i\varphi\left(\sum_{j=1}^n\omega(i,j)\lambda_jx_j\right)
\leq\sum_{i=1}^m\sum_{j=1}^n\mu_i\omega(i,j)\lambda_j\varphi(x_j)$$
$$
=\sum_{j=1}^n\left(\sum_{i=1}^m\mu_i\omega(i,j)\right)\lambda_j\varphi(x_j)
=\sum_{j=1}^n\lambda_j\varphi(x_j),
$$
and
$$
\sum_{i=1}^m\mu_i\varphi\left(\sum_{j=1}^n\omega(i,j)\lambda_jx_j\right)
\geq\varphi\left(\sum_{i=1}^m\sum_{j=1}^n\mu_i\omega(i,j)\lambda_jx_j\right)$$
$$
=\varphi\left(\sum_{j=1}^n\left(\sum_{i=1}^m\mu_i\omega(i,j)\right)\lambda_jx_j\right)
\nonumber\\
=\varphi\left(\sum_{j=1}^n\lambda_jx_j\right),$$
and (3) follows.\\
The inequalities in (4) follow from (3) by taking
$m=n,~\mu_i=\lambda_j=\frac{1}{n}~ \mbox{and}
~\omega(i,j)=na_{ij}~(i,j=1,\cdots,n)$.\\\\
Next, we give a refinement of discrete Jensen's inequality via
two weights functions.\\\\
{\bf Theorem 2.2.} {\it If $\omega_1$ and $\omega_2$ are two
weight functions, then
\begin {itemize} \item[\mbox{$(i)$}]
\begin{eqnarray}
\varphi\left(\sum_{j=1}^n\lambda_jx_j\right)\leq
\phi_{\omega_1,\omega_2}(t) \leq
\sum_{j=1}^n\lambda_j\varphi(x_j)~~~~~~(0\leq t \leq 1).
\end{eqnarray}
\item[\mbox{$(ii)$}] For each $i$, the function
$$
t\longrightarrow\varphi\left(\sum_{j=1}^n[(1-t)\omega_1(i,j)+t\omega_2(i,j)]\lambda_jx_j\right)
~~~~~~(0\leq t \leq 1),
$$
and so, $\phi_{\omega_1,\omega_2}$ is convex.\\
\item[\mbox{$(iii)$}]
\begin{eqnarray}
\varphi\left(\sum_{j=1}^n\lambda_jx_j\right)\leq
\int_{0}^1\phi_{\omega_1,\omega_2}(t)d t \leq
\sum_{j=1}^n\lambda_j\varphi(x_j).
\end{eqnarray}
In particular, if $C$ is an interval of ${\Bbb R}$,
\begin{eqnarray}
\varphi\left(\sum_{j=1}^n\lambda_jx_j\right)
\leq\sum_{i=1}^m\mu_iA
\left(\varphi;\sum_{j=1}^n\omega_1(i,j)\lambda_jx_j,\sum_{j=1}^n
\omega_2(i,j)\lambda_jx_j\right)\leq
\sum_{j=1}^n\lambda_j\varphi(x_j),
\end{eqnarray}
where the arithmetic mean $A$ is defined for an integrable
function $f$ over an interval with end points $a$ and $b$, by
\begin{eqnarray}
A(f; a, b)=\frac{1}{b-a}\int_{a}^b f(x)d x .
\end{eqnarray}
$($We set $A(f; a, a)=f(a)$$.)$
\item[\mbox{$(iv)$}] Let $p_i\geq 0$ with
$P_k=\sum _{i=1}^k p_i>0$, and $t_i$ be in $[0,1]$ for all
$i=1,2,\cdots,k$. Then
\begin{eqnarray}
\varphi\left(\sum_{j=1}^n\lambda_jx_j\right)\leq
\phi_{\omega_1,\omega_2}\left(\frac{1}{P_k}\sum_{i=1}^k p_it_i
\right) \leq\frac{1}{P_k}\sum_{i=1}^k
p_i\phi_{\omega_1,\omega_2}(t_i)\leq
\sum_{j=1}^n\lambda_j\varphi(x_j),
\end{eqnarray}
\end{itemize}
which is a discrete version of Hadamard's result.}\\
{\it Proof}
\begin{itemize}
\item[\mbox{$(i)$}] Since for each $t$ in $[0,1]$,
$$
(i,j)\longrightarrow(1-t)\omega_1(i,j)+t\omega_2(i,j)~~~~~~(1\leq
i\leq m,~1\leq j\leq n)
$$
is a weight function, (5) follows from (3).\\
\item[\mbox{$(ii)$}] Let $\alpha,~\beta\geq 0$ with $\alpha+\beta=1$ and
$t_1,t_2$ be in $[0,1]$. For each $i$ with $1\leq i\leq m$, we
have
\begin{eqnarray}
&~&\varphi\left(\sum_{j=1}^n[(1-\alpha t_1-\beta t_2
)\omega_1(i,j)+(\alpha t_1+\beta
t_2)\omega_2(i,j)]\lambda_jx_j\right)\nonumber\\
&=&
\varphi\left(\alpha\sum_{j=1}^n[(1-t_1)\omega_1(i,j)+t_1\omega_2(i,j)]\lambda_jx_j
+\beta\sum_{j=1}^n[(1-t_2)\omega_1(i,j)+t_2\omega_2(i,j)]\lambda_jx_j\right)\nonumber\\
&\leq&\alpha\varphi\left(\sum_{j=1}^n[(1-t_1)\omega_1(i,j)+t_1\omega_2(i,j)]\lambda_jx_j\right)
+\beta\varphi\left(\sum_{j=1}^n[(1-t_2)\omega_1(i,j)+t_2\omega_2(i,j)]\lambda_jx_j\right),
\nonumber
\end{eqnarray}
and (ii) follows.\\
\item[\mbox{$(iii)$}] $\phi_{\omega_1,\omega_2}$ being bounded and convex on
$[0,1]$ is Riemann
integrable on $[0,1]$, and so by integrating, (6) follows from (5).\\
In particular, if $C$ is an interval of ${\Bbb R}$, then by the
change of variables
$$
u=\sum_{j=1}^{n}[(1-t)\omega_1(i,j)+t\omega_2(i,j)]\lambda_jx_j,
$$
we have
\begin{eqnarray}
\int_0^1\phi_{\omega_1,\omega_2}(t)d
t&=&\sum_{i=1}^m\mu_i\int_0^1\varphi
\left(\sum_{j=1}^{n}[(1-t)\omega_1(i,j)+t\omega_2(i,j)]\lambda_jx_j\right)dt\nonumber\\
&=&\sum_{i=1}^m\mu_iA\left(\varphi;
\sum_{j=1}^{n}\omega_1(i,j)\lambda_jx_j,\sum_{j=1}^{n}\omega_2(i,j)\lambda_jx_j\right)
, \nonumber
\end{eqnarray}
which by substituting it in (6), we get (7).\\
\item[\mbox{$(iv)$}] The first and the third inequalities in (9) are obvious from (5),
and the second inequality follows from Jensen's inequality applied
for the convex function $\phi_{\omega_1,\omega_2}$.
\end{itemize}
{\bf Corollary 2.3.} {\it If $B=[b_{ij}]_{n\times n}$ and
$C=[c_{ij}]_{n\times n}$ are two double stochastic matrices, then
$[5]$,
\begin{eqnarray}
\varphi\left( \frac{x_1+\cdots+x_n}{n}\right)\leq\phi_{B,C}(t)\leq
\frac{\varphi(x_1)+\cdots+\varphi(x_n)}{n}~~~~~~(0\leq t \leq 1),
\end{eqnarray}and
\begin{eqnarray}
\varphi\left(
\frac{x_1+\cdots+x_n}{n}\right)\leq\int_{0}^1\phi_{B,C}(t)d t\leq
\frac{\varphi(x_1)+\cdots+\varphi(x_n)}{n}.
\end{eqnarray}
In particular, if $C$ is an interval of ${\Bbb R}$, then
\begin{eqnarray}
\varphi\left( \frac{x_1+\cdots+x_n}{n}\right)
\leq\frac{1}{n}\sum_{i=1}^n A\left(\varphi;
\sum_{j=1}^{n}b_{ij}x_j,\sum_{j=1}^{n}c_{ij}x_j\right)
\leq\frac{\varphi(x_1)+\cdots+\varphi(x_n)}{n},
\end{eqnarray}
where $A$ is defined by $(8)$.}\\
{\it Proof}. Take
$m=n,~\mu_i=\lambda_j=\frac{1}{n},~\omega_1(i,j)=nb_{ij}~\mbox{and}~\omega_2(i,j)=nc_{ij}
~(i,j=1,\cdots,n)$ in (5), (6) and (7).\\
\section{Applications}
Throughout this section, we use the terminologies and results of
the above sections, and as before, we suppose that $\omega_1$ and
$\omega_2$ are two weight functions, and $B=[b_{ij}]_{n\times n}$
and $C=[c_{ij}]_{n\times n}$ are two double stochastic matrices.
\\
\\{\bf Theorem 3.1.} {\it Let $x_1, x_2, \cdots, x_n$ be $n$
positive numbers. Then, we have
\begin{eqnarray}
\prod_{j=1}^n x_j^{\lambda_j}\leq\prod_{i=1}^m \left[I\left(
\sum_{j=1}^{n}\omega_1(i,j)\lambda_j
x_j,\sum_{j=1}^{n}\omega_2(i,j)\lambda_jx_j\right)\right]^{\mu_i}\leq
\sum_{j=1}^n \lambda_jx_j,
\end{eqnarray}
where the identric mean $I$ is defined for each $a,b>0$ by
\begin{eqnarray}
I(a,b)= \left\{
\begin{array}{cl}
a &{\rm if}~~a=b,\\
\frac{1}{e}\left(\frac{b^b}{a^a}\right)^\frac{1}{b-a} &{\rm
if}~~a\not=b.
\end{array}
\right.
\end{eqnarray}
In particular $[5]$,
\begin{eqnarray}
\sqrt[n]{x_1 x_2 \cdots x_n} \leq
\sqrt[n]{\prod_{i=1}^nI\left(\sum_{j=1}^n b_{ij}x_j,\sum_{j=1}^n
c_{ij}x_j\right)}\leq\frac{x_1+ x_2+\cdots+ x_n}{n}.
\end{eqnarray}}
{\it Proof}. The function
$\varphi:(0,\infty)\rightarrow\mathbb{R},~\varphi(x)=-\ln x$  is
convex and $A(\varphi; a, b)=-\ln I(a,b)~(a,b>0)$. So, we have
$$
\sum_{i=1}^m\mu_iA\left(\varphi;
\sum_{j=1}^n\omega_1(i,j)\lambda_jx_j,\sum_{j=1}^n\omega_2(i,j)\lambda_jx_j\right)=
-\ln\prod_{i=1}^m\left[ I\left(
\sum_{j=1}^n\omega_1(i,j)\lambda_jx_j,\sum_{j=1}^n\omega_2(i,j)\lambda_jx_j\right)\right]^
{\mu_i},
$$
which by substituting in (7) and taking into account that
$$
\varphi \left (\sum_{j=1}^n\lambda_jx_j \right )=-\ln\left
(\sum_{j=1}^n\lambda_jx_j \right )\hspace{1cm}
\mbox{and}\hspace{1cm}\sum_{j=1}^n\lambda_j\varphi(x_j)=-\ln
\left(\prod_{j=1}^n x_j^{\lambda_j}\right),
$$
we obtain (13).\\
The inequalities in (15) follow from (13) by taking
$m=n,~\mu_i=\lambda_j=\frac{1}{n},~\omega_1(i,j)=nb_{ij}~\mbox{and}~
\omega_2(i,j)=nc_{ij}~(i,j=1,\cdots,
n)$.\\
\\{\bf Theorem 3.2.} {\it If $x_j\in (0,1/2]~(j=1,\cdots,n)$, and
$A_n=\sum_{j=1}^n \lambda_jx_j$ and $G_n=\prod_{j=1}^n
x_j^{\lambda_j}$~$($also, $A'_n=\sum_{j=1}^n \lambda_j(1-x_j)$ and
$G'_n=\prod_{j=1}^n (1-x_j)^{\lambda_j}$~$)$ are the arithmetic
and geometric means of $x_1, \cdots, x_n$~$($of $1-x_1, \cdots,
1-x_n$~$)$ respectively, then we have the following refinement of
Ky Fan's inequality $\frac{A'_n}{A_n}\leq \frac{G'_n}{G_n}$ $[1]$:
\begin{eqnarray}
\frac{A'_n}{A_n}\leq\prod_{i=1}^{m}\left(
\frac{I\left(\sum_{j=1}^n\omega_1(i,j)\lambda_j(1-x_j),
\sum_{j=1}^n\omega_2(i,j)\lambda_j(1-x_j)\right)}
{I\left(\sum_{j=1}^n\omega_1(i,j)\lambda_jx_j,
\sum_{j=1}^n\omega_2(i,j)\lambda_jx_j\right)}\right)^
{\mu_i} \leq\frac{G'_n}{G_n},
\end{eqnarray}
where the Identric mean $I$ is defined as in $(14)$.\\
In particular $[6]$,
\begin{eqnarray}
\frac{A'_n}{A_n}\leq \sqrt[n]{\prod_{i=1}^{n}\left(
\frac{I\left(\sum_{j=1}^nb_{ij}(1-x_j),\sum_{j=1}^nc_{ij}(1-x_j)\right)}
{I\left(\sum_{j=1}^nb_{ij}x_j,\sum_{j=1}^nc_{ij}x_j\right)}\right)}
\leq\frac{G'_n}{G_n}.
\end{eqnarray}}
{\it Proof}. The function $\varphi(x)=\ln\frac{1-x}{x}$ is convex
on $(0,1/2]$, and $A(\varphi; a, b)=\ln
\frac{I(1-a,1-b)}{I(a,b)}~(0<a,b<1)$. So, we have
$$\sum_{i=1}^m\mu_iA\left(\varphi;
\sum_{j=1}^n\omega_1(i,j)\lambda_jx_j,\sum_{j=1}^n\omega_2(i,j)\lambda_jx_j\right)$$
$$=\ln \prod_{i=1}^{m}\left(
\frac{I\left(\sum_{j=1}^n\omega_1(i,j)\lambda_j(1-x_j),
\sum_{j=1}^n\omega_2(i,j)\lambda_j(1-x_j)\right)}
{I\left(\sum_{j=1}^n\omega_1(i,j)\lambda_jx_j,
\sum_{j=1}^n\omega_2(i,j)\lambda_jx_j\right)}\right)^{\mu_i},$$
which by substituting in (7) and taking into account that
$$
\varphi \left (\sum_{j=1}^n\lambda_jx_j \right )=\ln
\frac{A'_n}{A_n}\hspace{1cm}
\mbox{and}\hspace{1cm}\sum_{j=1}^n\lambda_j\varphi(x_j)=\ln
\frac{G'_n}{G_n},
$$
we get (16).\\
In particular, (17) follows from (16) by taking
$m=n,~\mu_i=\lambda_j=1/n,~\omega_1(i,j)=nb_{ij}~\mbox{and}~
\omega_2(i,j)=nc_{ij}~(i,j=1,\cdots,n)$.\\\\
{\bf Theorem 3.3.} {\it If $(X,\cal{A},\mu)$ is a measure space,
$p\geq 1$, and $f_1, f_2,\cdots,f_n$ belong to $L^p=L^p(\mu)$,
then we have
\begin{eqnarray}
\left\|\sum_{j=1}^n\lambda_jf_j\right\|_{p}^{p}
\leq\sum_{i=1}^{m}\mu_i
\left\|L_p^p\left(\sum_{j=1}^n\omega_1(i,j)\lambda_j|f_j|,
~\sum_{j=1}^n\omega_2(i,j)\lambda_j|f_j|\right)\right\|_1
\leq\sum_{j=1}^n\lambda_j\|f_j\|_p^p,
\end{eqnarray}
where the $p$-logarithmic mean is defined for $a,b\geq 0$, by
\begin{eqnarray} L_p(a,b)= \left\{
\begin{array}{cl}
a &{\rm if}~~a=b,\\
\left[\frac{b^{p+1}-a^{p+1}}{(p+1)(b-a)}\right]^{1/p}&{\rm
if}~~a\not=b.
\end{array}
\right.
\end{eqnarray}
In particular $[5]$,
\begin{eqnarray}
\left\|\frac{f_1+\cdots+f_n}{n}\right\|_{p}^{p}
\leq\frac{1}{n}\sum_{i=1}^{n}
\left\|L_p^p\left(\sum_{j=1}^{n}b_{ij}|f_j|,~\sum_{j=1}^{n}c_{ij}|f_j|\right)\right\|_1
\leq\frac{\|f_1\|_p^p+\cdots+\|f_n\|_p^p}{n},
\end{eqnarray}
and
\begin{eqnarray}
\left\|\frac{f_1+\cdots+f_n}{n}\right\|_{p}^{p}
\leq\frac{1}{n}\sum_{i=1}^{n}
\left\|L_p^p(|f_i|,~|f_{n+1-i}|)\right\|_1
\leq\frac{\|f_1\|_p^p+\cdots+\|f_n\|_p^p}{n}.
\end{eqnarray}}
{\it Proof}. We consider the convex function
$\varphi:L^p\rightarrow {\Bbb{R}},~\varphi(f)=\|f\|_p^p$. Clearly,
the function $X\times [0,1]\rightarrow {\Bbb{R}}$ with
$$
(x,t)\rightarrow\sum_{j=1}^n[(1-t)\omega_1(i,j)+t\omega_2(i,j)]\lambda_jf_j(x),
$$
is product-measurable. Since
$$
\left|\sum_{j=1}^n\lambda_jf_j\right|\leq\sum_{j=1}^n\lambda_j|f_j|,
$$
and the $L^p-$norms of $f_j$ and $|f_j|$ are equal~$(j=1,\cdots,
n)$, it is sufficient to assume $f_j\geq 0~(j=1,\cdots, n)$. Now,
using Fubini's theorem and applying the change of variables
$$u=\sum_{j=1}^n[(1-t)\omega_1(i,j)+t\omega_2(i,j)]\lambda_jf_j(x),$$
we have
\begin{eqnarray}
\int_{0}^{1}\phi_{\omega_1,\omega_2}(t)d
t&=&\sum_{i=1}^{m}\mu_i\int_{0}^{1}
\left\|\sum_{j=1}^{n}[(1-t)\omega_1(i,j)+t\omega_2(i,j)]\lambda_jf_j\right\|_p^pdt\nonumber\\
&=& \sum_{i=1}^{m}\mu_i\int_{0}^{1}
\int_X\left(\sum_{j=1}^{n}[(1-t)\omega_1(i,j)+t\omega_2(i,j)]\lambda_jf_j(x)\right)^pd\mu(x)dt
\nonumber
\end{eqnarray}
\begin{eqnarray}
&=&\sum_{i=1}^{m}\mu_i\int_X
\int_{0}^{1}\left(\sum_{j=1}^{n}[(1-t)\omega_1(i,j)+t\omega_2(i,j)]\lambda_jf_j(x)\right)^p
dtd\mu(x)
\nonumber\\
&=&\sum_{i=1}^{m}\mu_i\int_X
L_p^p\left(\sum_{j=1}^n\omega_1(i,j)\lambda_jf_j,~\sum_{j=1}^n\omega_2(i,j)\lambda_jf_j\right)
d\mu
\nonumber\\
&=&\sum_{i=1}^{m}\mu_i
\left\|L_p^p\left(\sum_{j=1}^n\omega_1(i,j)\lambda_jf_j,
~\sum_{j=1}^n\omega_2(i,j)\lambda_jf_j\right)\right\|_1, \nonumber
\end{eqnarray}
which by substituting in (6) with $f_j$ instead of $x_j$, it yields (18).\\
In particular, (20) follows from (18) by taking
$m=n,~\mu_i=\lambda_j=\frac{1}{n},~\omega_1(i,j)=nb_{ij},~\omega_2(i,j)=nc_{ij}
~(i,j=1,\cdots,n)$.\\
Finally, (21) follows from (20) by taking $b_{ij}=\delta_{ij}$ and
$c_{ij}=\delta_{i,n+1-j}~(i,j=1,\cdots,n)$, where $\delta_{ij}$ is
the Kronecker delta.\\\\
{\bf Corollary 3.4.} {\it If $x_1,x_2,\cdots,x_n$ are nonnegative
real numbers and $p\geq 1$, then
\begin{eqnarray}
\sum_{j=1}^n\lambda_j^px_j^p \leq\sum_{i=1}^{m}\sum_{j=1}^{n}\mu_i
L_p^p\left(\omega_1(i,j)\lambda_jx_j,
~\omega_2(i,j)\lambda_jx_j\right) \leq\sum_{j=1}^n\lambda_jx_j^p.
\end{eqnarray}
As a consequence, if $p$ is a positive integer, then
\begin{eqnarray}
n^{2-p}
\leq\frac{1}{p+1}\sum_{i,j=1}^{n}\sum_{k=0}^pb_{ij}^kc_{ij}^{p-k}
\leq n,
\end{eqnarray}
and
\begin{eqnarray}
n^{2-p}
\leq\frac{1}{p+1}\left(\sum_{i,j=1}^{n}b_{ij}^p+\sum_{i=1}^{n}\sum_{k=0}^{p-1}b_{ii}^k
\right)\leq n.
\end{eqnarray}}
{\it Proof}. The inequalities in (22) follow from (18) by taking
$X=\{1,2,\cdots,n\}$ with counting measure and
$f_j=x_j\chi_{\{j\}}$, where $\chi_{\{j\}}$ is the characteristic function of
$\{j\}~(j=1,2,\cdots,n)$.\\
Now, (23) follows from (22) by taking
$m=n,~\mu_i=\lambda_j=\frac{1}{n},~\omega_1(i,j)=nb_{ij},~\omega_2(i,j)=nc_{ij}$,
$x_j=1~(i,j=1,\cdots,n)$ and expanding the $p$-logarithmic mean.\\
Finally, (24) follows from (23) by taking
$c_{ij}=\delta_{ij}~(i,j=1,\cdots,n)$.
\\\\{\bf Remark 3.5.} (i) Let ${f_n}$ be a sequence in $L^p(\mu)~(p\geq1)$ converging
with the $L^p$-norm and point-wise to an element $f$ of
$L^p(\mu)$. Then, using Fatu's lemma and Cesaro's summability
theorem, we have
$$
\|f\|_p^p\leq
\liminf_{n\rightarrow\infty}\left\|\frac{f_1+\cdots+f_n}{n}\right\|_{p}^{p}
\leq\lim_{n\rightarrow\infty}\frac{\|f_1\|_p^p+\cdots+\|f_n\|_p^p}{n}=\|f\|_p^p,
$$
and so by (21),
\begin{eqnarray}\lim_{n\rightarrow\infty}\frac{1}{n}\sum_{i=1}^{n}
\|L_p^p(|f_i|,~|f_{n+1-i}|)\|_1=\|f\|_p^p~~~~~~(p\geq1).
\end{eqnarray}
(ii) Let $(X,\cal{A},\mu)$ be a finite measure space and ${\cal
M}$ be the vector space of all measurable functions on $X$ with
point-wise operations [2]. The set $C$, consisting of all
nonnegative measurable functions on $X$, is a convex subset of
${\cal M}$. Since the function $t\rightarrow\frac{t}{1+t}
~~(t\geq0)$ is concave, the mapping
$\varphi:C\rightarrow\mathbb{R}$ with
\begin{eqnarray}
\varphi(f)=\int_X\frac{f}{1+f}d \mu\hspace{1cm}(f\in C)
\end{eqnarray}
is concave.\\\\
{\bf Theorem 3.6.} {\it With the notations of $(ii)$ of Remark
3.5, if $f_1,\cdots,f_n$ belong to $C$ and $\varphi$ is as in
$(26)$, then
\begin{eqnarray}
&&\sum_{j=1}^{n}\lambda_j\varphi(f_j)\\
&\leq&\mu(X)-\sum_{i=1}^{m}\mu_i \left\|L^{-1}\left(
1+\sum_{j=1}^n\omega_1(i,j)\lambda_jf_j,
~1+\sum_{j=1}^n\omega_2(i,j)\lambda_jf_j\right)\right\|_1 \nonumber\\
&\leq&\varphi\left(\sum_{j=1}^{n}\lambda_jf_j\right),\nonumber
\end{eqnarray}
where the logarithmic mean $L$ is defined for each $a, b
>0$, by
\begin{eqnarray}
L(a,b)= \left\{
\begin{array}{cl}
a &{\rm if}~~a=b,\\
\frac{b-a}{\ln b-\ln a}&{\rm if}~~a\not=b.
\end{array}
\right.
\end{eqnarray}
In particular $[5]$,
\begin{eqnarray}
&&\frac{\varphi(f_1)+\cdots+\varphi(f_n)}{n}\\
&\leq&\mu(X)-\frac{1}{n}\sum_{i=1}^{n}\left\|L^{-1}\left(
1+\sum_{j=1}^nb_{ij}f_j, ~1+\sum_{j=1}^nc_{ij}f_j\right)\right
\|_1
\nonumber\\
&\leq&\varphi\left(\frac{f_1+\cdots+f_n}{n}\right) ,\nonumber
\end{eqnarray}}
{\it Proof}. Clearly, the mapping
$$
(x,t)\rightarrow\sum_{j=1}^n[(1-t)\omega_1(i,j)+t\omega_2(i,j)]\lambda_jf_j(x)
$$
on $X\times [0,1]$ is product-measurable. Since $\varphi$ is
concave, $-\varphi$ is convex, and so by (6), we have
\begin{eqnarray}
&&\sum_{j=1}^{n}\lambda_j\varphi(f_j)\leq
\int_0^1\phi_{\omega_1,\omega_2}(t)d t
\leq\varphi\left(\sum_{j=1}^{n}\lambda_jf_j\right).
\end{eqnarray}
But
\begin{eqnarray}
&~&\int_{0}^{1}\phi_{\omega_1,\omega_2} (t)d t\nonumber\\
&=& \sum_{i=1}^{m}\mu_i\int_{0}^{1}\int_X\frac
{\sum_{j=1}^n[(1-t)\omega_1(i,j)+t\omega_2(i,j)]\lambda_jf_j(x)}
{1+\sum_{j=1}^n[(1-t)\omega_1(i,j)+t\omega_2(i,j)]\lambda_jf_j(x)}d\mu(x)dt\nonumber\\
&=& \sum_{i=1}^{m}\mu_i\int_X\int_0^1\frac
{\sum_{j=1}^n[(1-t)\omega_1(i,j)+t\omega_2(i,j)]\lambda_jf_j(x)}
{1+\sum_{j=1}^n[(1-t)\omega_1(i,j)+t\omega_2(i,j)]\lambda_jf_j(x)}dtd\mu(x)\nonumber
\end{eqnarray}
\begin{eqnarray}&=&\sum_{i=1}^{m}\mu_i\int_X\frac{1}{\sum_{j=1}^n[\omega_2(i,j)-\omega_1(i,j)]\lambda_jf_j(x)}
\int_{\sum_{j=1}^n\omega_1(i,j)\lambda_jf_j(x)}^
{\sum_{j=1}^n\omega_2(i,j)\lambda_jf_j(x)}(1-\frac{1}{1+t})dtd\mu(x)\nonumber\\
&=&\mu(X)-\sum_{i=1}^{m}\mu_i\int_X
\frac{1}{\sum_{j=1}^n[\omega_2(i,j)-\omega_1(i,j)]\lambda_jf_j}
\ln\frac{1+\sum_{j=1}^{n}\omega_2(i,j)\lambda_jf_j}
{1+\sum_{j=1}^{n}\omega_1(i,j)\lambda_jf_j}d \mu\nonumber\\
&=&\mu(X)-\sum_{i=1}^{m}\mu_i
\left\|L^{-1}\left(1+\sum_{j=1}^n\omega_1(i,j)\lambda_jf_j,
~1+\sum_{j=1}^n\omega_2(i,j)\lambda_jf_j\right)\right\|_1
,\nonumber
\end{eqnarray}
which by substituting in (30), we obtain (27).\\
The inequalities in (29) follow from (27) by taking
$m=n,~\mu_i=\lambda_j=\frac{1}{n},~\omega_1(i,j)=nb_{ij}~\mbox{and}~\omega_2(i,j)=nc_{ij}
~(i,j=1,\cdots,n)$.\\
\begin{center}
{REFERENCES}
\end{center}
\begin{enumerate}
\item E. F. Beckenbach and R. Bellman, {\it Inequalities}, Springer-Verlag, Berlin, 1961.
\item S. Berberian, {\it Lectures in functional analysis and operator theory},
Springer, New York-Heidelberg-Berlin, 1974.
\item S.S. Dragomir, On some refinement of Jensen's inequality and applications,
{\it Utilitas Mathematica}, {\bf 43}(1993), 235-243.
\item S.S. Dragomir, J.E. Pe\v{c}ari\'{c} and L.E. Persson,
Properties of some functionals related to Jensen's inequality,
{\it Acta Math. Hungar.}, {\bf 70} (1-2)(1996), 129-143.
\item J. Rooin, Some aspects of convex
functions and their applications, {\it Jipam}, Vol 2, Issue 1,
Article 4 (2001).
\item J. Rooin, Some refinements of
Ky Fan's and Sandor's inequalities, {\it Southeast Asian Bulletin
of Mathematics} (in press).
\end{enumerate}
\end{document}